\documentclass[12pt]{article}
\usepackage{epsfig,amsfonts,amsmath}
\def\Bbb R{{\rm \bf R}}
\def\proclaim#1{\vskip2mm{\bf #1}\em}
\def\endproclaim{\em \vskip2mm}
\def\tag#1{\eqno(#1)}
\def\gathered{\begin{array}{c}}
\def\endgathered{\end{array}}
\def\text{\mbox}

\begin{document}

\title {Trusted frequency region of convergence for the enclosure method
in an inverse heat equation }
\author{Masaru Ikehata\footnote{
Laboratory of Mathematics,
Institute of Engineering,
Hiroshima University,
Higashi-Hiroshima 739-8527, Japan.  e-mail address: ikehata@amath.hiroshima-u.ac.jp}
and
Kiwoon Kwon\footnote{Corresponding author: Department of Mathematics, Dongguk University-Seoul, 100715 Seoul, South Korea.
e-mail address:
kwkwon@dongguk.edu}
}
\maketitle

\begin{abstract}
This paper is concerned with the numerical implementation
of a formula in the enclosure method as applied to a prototype inverse initial
boundary value problem for thermal imaging in a one-space dimension.
A precise error estimate of the formula is given
and the effect on the discretization of the used integral
of the measured data in the formula is studied.  
The formula requires a large frequency to converge; however, 
the number of time interval divisions grows exponetially as the frequency increases.
Therefore, for a given number of divisions, we fixed the trusted frequency region of convergence
with some given error bound.  
The trusted frequency region is computed theoretically using  theorems provided in 
this paper and is numerically implemented for various cases. 
\noindent
AMS: 35R30

\end{abstract}


\section{Introduction}


Thermal imaging is described as follows: given a heat flux on the surface of an object and a measured surface temperature, determine the internal thermal properties of the object 
or the shape of some unknown inaccessible portion of the boundary \cite{BC}.
$$
\left\{
\begin{array}{ll}
\displaystyle
u_t=\triangle u & \text{in}\,\Omega\times {\mathbb R^+},
\\
\\
\displaystyle
\frac{\partial u}{\partial n}=f & \text{on}\,\partial\Omega\times{\mathbb R^+},
\\
\\
\displaystyle
u(x,0) = u_0(x) & \text{on}\,\Omega\times\{t=0\}.
\end{array}
\right.
\tag {1.1}
$$

Let the support of $f$ and the measurement set be contained in the known boundary 
$\Gamma\subset \partial\Omega$ . Thermal imaging is redescribed as determining the part of $\partial \Omega$ such that $ f=0 $, which means a perfectly insulating boundary. The problem is applied to identify back surface corrosion and damage, such as the use of infrared thermography to find burn injuries and the selection of donor sites for skin grafts.

It is reported in \cite{BC} that if
$$
\left\{
\begin{array}{l}
\displaystyle
\Omega_1=[0,2\pi]\times[0,\pi],
\\
\\
\displaystyle
\Omega_2=\Omega_1 \setminus \left[\frac {2\pi}{3}, \frac{4\pi}3\right]\times\left[0,\frac{2\pi}{3}\right],
\\
\\
\displaystyle
u(t,x,y)=e^{-\frac 9 2 t} \cos\left(\frac 3 2 x\right) \cos\left(\frac 3 2 y\right),
\\
\\
\displaystyle
\Gamma =[0,2\pi]\times [y=\pi],
\end{array}
\right.
$$
$\Omega_1\setminus\Gamma$ and $\Omega_2\setminus\Gamma$ are two different 
unknown boundaries on which the Neumann data $f$ vanishes. This is an example of the nonuniqueness of the thermal imaging problem.

On the other hand, two uniqueness results are also reported in \cite{BC}. 
\begin{itemize}
\item{ If $u_0$ is constant and $u_1=u_2$ on $\Gamma\times(0,T)$, then we have $\Omega_1 = \Omega_2$ and $u_1=u_2$.}
\item{ If $u_0$ is nonconstant, special conditions are required for the uniqueness of $\Omega$ and $u$. That is, if  
$$
\Vert f(t,\cdot)\Vert_{L^2(\partial\Omega)}, 
\Vert f_t(t,\cdot)\Vert_{L^2(\partial\Omega)} < \infty, 
\int_{\partial\Omega} f(t,x) dS_x >G_0 >0 ,$$ 
then $u_1=u_2$ on $\Gamma\times(0,\infty)$ implies $\Omega_1 = \Omega_2$ and $u_1=u_2$. }
\end{itemize}

A boundary element method is presented for a linearised inverse problem of  (1.1), as a numerical method \cite{BC96}. On the other hand, in this paper, the enclousre method is used for the nonlinear inverse problem of (1.1).

In \cite{I4}, some inverse problems for the heat and wave equations were included in a one-space dimension, and 
the first author introduced the {\it enclosure method in a time domain}.
The enclosure method is an analytical method which has its origins in \cite{I1} and \cite{IINF}.
Therein, the governing equations are elliptic equations and the observation data are given by
a single set of Cauchy data and the Dirichlet-to-Neumann map, respectively.
The enclosure method developed in \cite{I4} can be considered as an extension of the concept in \cite{I1}
to include inverse problems in the time domain. 
See also \cite{IE0,IE2, IK1,IK2, IK3}.

It is reported that the numerical implementation of thermal imaging without any linearisation as in \cite{BC96} even in one -space dimensional case is not trivial \cite{KS}.  Let us consider the following one-dimensional thermal imaging problem with constant initial data. 

Let $0<a<\infty$ and $0<T<\infty$.  Given $f\in L^2(0,\,T)$ let $u$ be a solution
of the problem:
$$\left\{
\begin{array}{ll}
\displaystyle
u_t=u_{xx} &  \text{in}\,\,]0,\,a[\times\,]0,\,T[, \\
\\
\displaystyle
u_x(a,t)=0  &  \text{for}\,t\in\,]0,\,T[, \\
\\
\displaystyle
u_x(0,t)=f(t)  & \text{for}\,t\in\,]0,\,T[, \\
\\
\displaystyle
u(x,0)=0  & \text{in}\,]0,\,a[.
\end{array}
\right.
\tag {1.2}
$$
Note that, because the initial data is constant, we can choose any nonzero Neumann data $f$ for the uniqueness of the unknown perfect conducting boundary $a$: However, we impose some weak condition (1.4) for $f$ for the enclosure method to be valid. The solution class is the same as that in \cite{I4, IK1} which was obtained from \cite{DL}.

Let $\tau>0$ and 
$$\displaystyle
v(x,t)=e^{-\tau^2 t}e^{-x\tau}.
$$
This $v$ satisfies the backward heat equation $v_t+v_{xx}=0$ in $]0,\,\infty[\times\,]0,\,T[$.

The so-called {\it indicator function} for the enclosure method here takes the form
$$\begin{array}{ccl}
\displaystyle
I(\tau)&=&\int_0^T\left( u_x(0,t) v(0,t) - v_x(0,t)u(0,t)\right)dt
\\
\\
\displaystyle
&=& \tau \hat u (0,\tau) + \hat f (\tau)
\end{array}
\tag {1.3}
$$
where $u$ satisfies (1.2) and $\hat w (\tau)= \int_0^T e^{-\tau^2 t} w(t) dt$ which is a modified Laplace transform
with finite time interval $T$. Let us consider $\tau$ to be the frequency corresponding to the enclosure method.

Assume that there exist positive numbers $C_\mu$, $\mu$, and $\tau_0$ such that
$$
\displaystyle
\left\vert \hat f(\tau)\right\vert\ge C_\mu\,\tau^{-\mu},  \,\,\forall\tau\ge\tau_0.
\tag {1.4}
$$
Then, by \cite{I4}, we have the formula
$$\displaystyle
\lim_{\tau\longrightarrow\infty}
\frac{1}{-2\tau}\log\vert I(\tau)\vert
=a.
\tag {1.5}
$$

Note that (1.4) is a restriction of the strength of the heat flux at $t=0$ from below.
In particular, $f(t)$ cannot be $0$ at $t=0$ with infinite order.  
It is easy to see that 
condition (1.4) is satisfied
if $f\in L^2(0,\,T)$ satisfies one of the following conditions for some $\delta\in\,]0,\,T[$:

$\bullet$  $\exists C>0$ such that $f(t)\ge C$ a.e. in $]0,\,\delta[$.

$\bullet$  $f\in C^{1}([0,\,\delta])$
and $f(0)\not=0$;

$\bullet$  $f\in C^{l+1}([0,\,\delta])$ with $l\ge 1$
and $f^{(s)}(0)=0$ for all $s=0,1,\cdots,l-1$ and $f^{(l)}(0)\not=0$.


When $f(t)=t^r, r=0,1,2,\cdots$ and $\tau>\frac 1 {\sqrt{T}}$, we have
$$ 0< C_\mu \tau^\mu \le \hat f(\tau) \le r! \tau^\mu
	\tag{1.6}$$
where
$$
\mu = 2(r+1), \;\; C_\mu = r!\left(   1- \frac{\sum_{k=0}^r \frac 1 {k!}}{e} \right). 
$$

(1.5) extracts $a$ from $u(0,t)$ given at a.e. $t\in\,]0,\,T[$ for a fixed known $f$.
A naive extraction procedure of $a$ is: just fix a large $\tau$ and compute an approximation
of $a$ such as 
$$\displaystyle
\frac{\log\vert I(\tau)\vert}{-2}
\approx a\tau.
$$
by finding a linear function fitting some values of $\frac{\log\vert I(\tau)\vert}{-2}$ at $\tau=\tau_1,\cdots\,\tau_m$
in the least-square sense and compute its slope which will be a candidate for the approximation of $a$.
This idea has been introduced in \cite{IS} for the enclosure method \cite{IINF}
and tested using an analytical solution of the direct problem. See also \cite{IO} for the enclosure method \cite{I1}.
Therein a similar numerical method has been tested using a solution of the direct problem constructed by finite element method. However, in this paper, rather than using linear approximation, a direct computation will be used with precise error analysis.

In this paper, instead of using (1.5) we develop another formula which is mathemaitcally equivalent . That is,
$$
\displaystyle
\lim_{\tau\longrightarrow\infty}
a(\tau)
=a, \;\;
\tag {1.7}
$$
where
$$\displaystyle
a(\tau)=\frac{1}{-2\tau}\log\left\vert\frac{I(\tau)}{-2\hat{f}(\tau)}\right\vert.
$$
Note that, since $f\in L^2(0,T)$ and satisfies (1.4), we have 
$$\displaystyle
\lim_{\tau\longrightarrow\infty}
\frac{\log\vert\hat{f}(\tau)\vert}{\tau}=0.
$$
Therefore, (1.5) and (1.7) are mathematically equivalent for $f\in L^2(0,T)$ satisfying (1.4). 

Then, what is the advantage of using (1.7) rahter than (1.5)? The reason is the following asymptotic formula
as $\tau\longrightarrow\infty$:
$$\displaystyle
I(\tau)=-2\hat{f}(\tau)e^{-2\tau a}(1+O(\tau^{-1})),
\tag {1.8}
$$
$$\displaystyle
a(\tau)-a 
= O(\tau^{-1}). 
\tag {1.9}
$$
Although the asymptotic convergence (1.8) is covered in \cite{I4} for equations that are more general than (1.2), the formula 
(1.5), instead of (1.7), is used for the numerical approximation; such inconsistant use of a formula makes the numerical scheme have not optimal order of convergence, even if a direct method is used. In this paper, we reprove (1.8),
prove the approximation error (1.9), and derive a numerical scheme based on (1.7). That is, we introduce a numerical method based on (1.7) instead of (1.5). This approach would enable us to perform error analysis indicating the convergence order depending on the frequency $\tau$, final time $T$, and the Neumann data $f$, which would not be given when we use (1.5).  This is the main reason for constructing the present numerical method based on (1.7), instead of on (1.5). In detail, we could have the following theorem:

\proclaim{\noindent Theorem 1.1.}
Assume that we know two positive constants $a_L$ and $a_U$ such that
$$\displaystyle
a_L\le a\le a_U.
$$
Assmume that $f\in L^\infty[0,\,T]$. Further, assume that there exists a positive number $\tau_0$ such that (1.4) holds for all $\tau\ge\tau_0$, 
$$\displaystyle
\tau_0\ge
\frac{3a_U}{4T}\left(1+\sqrt{1+\frac{8T\mu}{9 a_U^2}}\right),
\tag {1.10}
$$
and 
$$\displaystyle
\frac{e^{-T\tau_0^2+3a_U\tau_0}\tau_0^\mu}{2C_\mu} C_T \Vert f\Vert_{L^\infty[0,T]}
\le \epsilon <1,
\tag {1.11}
$$
where $C_T$ is given in (2.10).
Then, for all $\tau\ge\tau_0$ we have
$$\begin{array}{c}
\displaystyle
\vert a(\tau)-a\vert
\le
\frac{-\log(1-e^{-2a_L\tau})}{2\tau}
+\frac{ C_T \Vert f\Vert_{L^\infty[0,T]}}{4C_\mu(1-\epsilon)} 
e^{-T\tau^2+3a_U\tau}\tau^{\mu-1}  \\
\\
\displaystyle
\le
\frac{-\log(1-e^{-2a_L\tau_0})}{2\tau_0}
+\frac{\epsilon}{2\tau_0(1-\epsilon)}.
\end{array}
\tag {1.12}
$$
\endproclaim

Conditions (1.4), (1.10), and (1.11) are the criteria for the choice of $\tau_0$ when $a_U,a_L,C_\mu,\mu,C_T,\Vert f\Vert$  are known. This result ensures the accuracy of the approximation $a(\tau)$ exactly for $a$ for all $\tau\ge\tau_0$. Thus, the problem becomes that of how to compute $a(\tau)$ as precisely as possible from observation data.

In the computation of $a(\tau)$ in (1.3) and (1.7), we need $u(0,t)$ for all $t\in]0,\,T[$. However, in practice, it is not possible to know $u(0,t)$ for all $t\in]0,\,T[$. Here, we consider how to compute $a(\tau)$ approximately from temperatures $u(t_j), t_j=\frac{jT}{N_t}, j=0,\cdots, N_t$ equidistantly sampled at $N_t$ discrete times taken from time interval $[0,\,T]$.

Let 
$$\displaystyle 
Q_L\left(\int_0^T g(t)dt\right) = \frac 1 L \sum_{k=1}^{L-1} g\left(\frac{kT}{L}\right) + \frac {g(0)+g(T)}{2L}
$$
 denote the trapezoidal rule for the integral of a continuous function $g$ over $[0,\,T]$ with
$L$ equidistant subdivision.  It is well known (see \cite{A}) that if $g$ is twice continuously differentiable, then the error has the estimate
$$\displaystyle
\left\vert\int_0^T g(t)dt-Q_L\left(\int_0^T g(t)dt\right)\right\vert
\le
\frac{T^3}{12L^2}\sup_{0\le t\le T}\vert g''(t)\vert.
\tag {1.13}
$$

Therefore, another issue that would have to be considered for the numerical implementation of (1.7) is the effect of the division number $N_t$ for the time interval $[0,T]$.
When the trapezoidal rule is used for $\hat u(\tau)$, it becomes possible to define the following:
\[
\left\{
\begin{array}{l}
\displaystyle
Q_{N_t}(\hat u(\tau))= 
\frac 1{N_t}\sum_{j=1}^{N_t-1} e^{ \frac {-jT\tau^2}{N_t}} u(0,\frac {jT}{N_t}) + \frac{u(0,0)+u(0,T)}{2 N_t},
\\
\\
\displaystyle
I_{N_t}(\tau)=\tau Q_{N_t}(\hat u(\tau)) +\hat f(\tau),
\\
\\
\displaystyle
a_{N_t}(\tau)=\frac 1{-2\tau} \log\left\vert \frac{I_{N_t}(\tau)}{-2\hat f(\tau)}\right\vert.
\end{array}
\right.
\tag{1.14}
\]

As the error (1.13) of the trapezoidal rule $Q_{N_t}(\hat u(\tau))$ depends on  $\frac 1 {N_t^2}$ and $\tau^4$, because of the second derivative of $e^{-\tau^2 t}u(0,t)$, the resulting error between $a$ and the approximation $a_{N_t}(\tau)$ is proportional to  $N_t^{-2}$ and $e^{2a\tau}\tau^4$. Therefore, for the approximation $a_{N_t}(\tau)$ to converge to $a$, it is required that $N_t$ is proportional to $e^{a \tau} \tau^{\frac{5+\mu}2+\delta}$ for a relatively large $\tau$ with some positive $\delta$ by the following Theorem 1.2, resulting in a numerically very expensive method. Remind that the norm for  the Sobolev space $W^{2,\infty}[0,T]$ is defined by  
$$
 \Vert f\Vert_{W^{2,\infty}[0,T]} = \max_{0\le s\le T} ( \vert f(s) \vert, \vert f_t(s) \vert, \vert f_{tt}(s) \vert ).
$$
\proclaim{\noindent Theorem 1.2.}
Under the conditions of Theorem 1.1, we further assume that $f\in W^{2,\infty}[0,\,T], \; f(0)=f'(0)=0$ and   
$\tau_0$ also satisfies 
$$\displaystyle
\frac{T^3C_{\max}\Vert f\Vert_{W^{2,\infty}[0,T]}}{24 C_\mu (1-\epsilon)}\left(\frac{1}{\tau_0^{\delta}}+\frac{1}{\tau_0^{2+\delta}}\right)^2\le\eta <1,
\tag {1.15}
$$
where $C_{max}$ is given in (2.3).
Then, it holds that, for all $\tau\ge\tau_0$ and $N_t\ge N_t^{\delta}(\tau):=[e^{a_U\tau}\tau^{(5+\mu+2\delta)/2}]+1$
$$\displaystyle
\vert a_{N_t}(\tau)-a(\tau)\vert
\le
\frac{T^3C_{\max}\Vert f\Vert_{W^{2,\infty}[0,T]}}{48 C_\mu (1-\eta)(1-\epsilon)}\frac 1\tau\left(\frac{1}{\tau^{\delta}}+\frac{1}{\tau^{2+\delta}}\right)^2
\le
\frac{\eta}{2\tau_0(1-\eta)}.
\tag {1.16}
$$
\endproclaim

Summing up the assumptions of Theorem 1.1 and Theorem 1.2, we should choose  $f\in W^{2,\infty}[0,\,T]$ satifying $f(0)=f'(0)=0$ and (1.4), and $\tau_0=\tau_0(\epsilon,\eta,\delta)$ satisfying (1.10),(1.11) and (1.15). If $N_t$ is greater than $N_t^\delta(\tau_0)$, we could define
\[
     \tau_{max} = \tau_{max}(N_t) := {\rm argmax} \{\tau\ge \tau_0 | N_t^\delta(\tau) \le N_t\} 
\]

Let us considerl $[\tau_0, \tau_{max}]$, the trusted frequency region with the error bound given in the following 
Theorem 1.3, which is simply derived from Theorems 1.1 and 1.2. 

\proclaim{\noindent Theorem 1.3}
Let the assumptions of Theorems 1.1 and 1.2 hold. Then, there exist $\tau_0$ 
and $\tau_{max}$ 
such that
\[
|a-a_{N_t}(\tau)| \le -\frac{\log(1+e^{-2 a_L\tau})}{2\tau} + \frac\epsilon{2\tau(1-\epsilon)}
+\frac\eta{2\tau(1-\eta)}
\]
for $\tau\in [\tau_0,\tau_{max}]$.
\endproclaim


In Section 2, Lemmas will be stated and proved before Theorems 1.1 and 1.2 are proved in Section 3. The numerical implementation is presented in Section 4.

\section{Lemmas}

The Riemann-Zeta function is defined as follows: 
$$ 
              \zeta(r) = \sum_{n=1}^\infty n^{-r}.
$$
It is well-known that  $\zeta(r),r=2,3,4,\cdots$ is a bounded real number and for even number $r=2k, k=1,2,\cdots$  
$$       \zeta(2k) = \frac{(-1)^r B_{2k} (2\pi)^{2k}}{2(2n)!},  $$
where $B_{2k}$ is a Bernoulli number.
For example, we have the values :
$$ \zeta(2) = \frac{\pi^2}6, \; \zeta(4) = \frac{\pi^4}{90}, \; \zeta(6)=\frac{\pi^6}{945}. $$

\proclaim{\noindent Lemma 2.1}
The solution $u$ of (1.2) is represented by 
$$\begin{array}{c}
\displaystyle
u(x,t)=-\frac{a}3 f(t)\left(\frac 3{2a^2} x^2 - \frac 3 a x\right)
-\frac{1}{a}
\int_0^t f(t)dt\\
\\
\displaystyle
-\frac 2 a\sum_{k=1}^{\infty}\cos\left(\frac{k\pi}a x\right)\int_0^t e^{-\left(\frac{k\pi}{a}\right)^2(t-s)}f(s)ds.
\end{array}
\tag{2.1}
$$
\endproclaim
Inserting $x=0$, we have the following Dirichlet data:
$$
u(0,t)= -\frac{1}{a}\int_0^t f(t)dt
\displaystyle
-\frac 2 a\sum_{k=1}^{\infty}\int_0^t e^{-\left(\frac{k\pi}{a}\right)^2(t-s)}f(s)ds.
\tag{2.2}
$$
Note that if $f(t)\ge 0$, then $u(0,t)\ge 0$ also.\\

{\bf Proof of Lemma 2.1}\\
For the problem (1,1), $v$ and eigenpairs $(\lambda_k,\Psi_k),k=1,2,\cdots$ in 
Lemma 3.2 in \cite{BC} are as follows:  
\noindent
$$\displaystyle
v(x,t) = -\frac{f(t)}{6a}(3x^2 - 6ax + 2a^2), \,\,
\lambda_k=\left(\frac{k\pi}{a}\right)^2.\;\;
\Psi_k(x)=\sqrt{\frac{2}{a}}\cos\,\frac{k\pi}{a}x.
$$

A direct computation yields
$$\left\{
\begin{array}{c}
\displaystyle
\int_0^a v(x,t)\Psi_k(x) dx =-f(t)\sqrt{\frac{2}{a}}\lambda_k^{-1},\\
\\
\displaystyle
\int_0^av_t(x,t)\Psi_k(x)dx = -f'(t)\sqrt{\frac{2}{a}}\lambda_k^{-1}.
\end{array}
\right.
$$

Using these computational results and Lemma 3.2 in \cite{BC}, we have the following representation formula: 
$$\begin{array}{c}
\displaystyle
u(x,t)=-\frac{a}3 f(t)\left(\frac 3{2a^2} x^2 - \frac 3 a x + 1\right)
-\frac{1}{a}
\int_0^t f(t)dt\\
\\
\displaystyle
+\frac 2 a\sum_{k=1}^{\infty}\cos\left(\frac{k\pi}a x\right)\frac 1 {\lambda_k} 
\left[e^{-\lambda_k t} f(0) -  \int_0^t e^{-\lambda_k(t-s)}f'(s)ds\right].
\end{array}
$$
Using integration by parts for the last integral and using $\zeta(2) = \frac {\pi^2}6$, we obtain  
equation (2.1).
$\square$\\

\proclaim{\noindent Lemma 2.2}
Assume that $f\in W^{2,\infty}[0,\,T]$ and $f(0)=f'(0)=0$. Then, we have
$$
    \Vert u(0,\cdot) \Vert_{W^{2,\infty}[0,T]}  \le C_{\max} \Vert f\Vert_{W^{2,\infty}[0,T]},
\tag{2.3}
$$
where 
$$
   C_{\max} =  \frac{\max(T,1)}{a_L} + \frac {a_U}3.
$$
If $f(t)=t^r,r=2,3,\cdots$, we have
$$
 \Vert u(0,\cdot) \Vert_{W^{2,\infty}[0,T]}   \le C_{\max,r},
\tag{2.4}
$$
where
$$\begin{array}{c}
\displaystyle
C_{\max,r} = \frac{T^{r-1}\max\left(\frac{T^2}{r+1}, T, r\right) }{a_L}
+ \frac{ a_U T^{r-2} \max\left(T^2,rT,r(r-1)\right)}  3.
\end{array}
$$
Further, if $T\ge r+1$, then
$$\displaystyle
C_{\max,r} = \frac{T^{r+1}}{a_L} + \frac{T^r a_U}3.
$$
\endproclaim

{\bf Proof of Lemma 2.2}\\
Using 
$$\displaystyle
\int_0^t e^{-\lambda_k(t-s)}\vert f(s)\vert ds
\le \Vert f\Vert_{L^\infty[0,t]} \int_0^te^{-2\lambda_k(t-s)}ds
\le \frac{\Vert f\Vert_{L^\infty[0,\,t]}}{\lambda_k}.
\tag{2.5}
$$
and $\zeta(2) = \frac {\pi^2} 6$, the upper bound of (2.2) is given by
$$\begin{array}{c}
\displaystyle
\vert u(0,t)\vert
\le
\frac 1 {a_L} \Vert f\Vert_{L^1[0,t]} + \frac {a_U} 3 \Vert f\Vert_{L^\infty[0,t]},
\end{array}
\tag{2.6}
$$
To enable a more convenient differntiation of $u(0,t)$ in (2.2), let us change (2.2) as follows by changing  $\eta = t-s$ in the last integral:
$$
u(0,t)= -\frac{1}{a}\int_0^t f(t)dt
\displaystyle
-\frac 2 a\sum_{k=1}^{\infty}\int_0^t e^{-\lambda_k\eta}f(t-\eta)d\eta.
\tag{2.7}
$$
By differentiating both sides of (2.7), we have
$$\begin{array}{l}
\displaystyle
u_t(0,t)=-\frac 1 a f(t)  - \frac 2 a \sum_{k=1}^\infty \left[ e^{-\lambda_k t}f(0) + \int_0^t e^{-\lambda_k\eta}f'(t-\eta)d\eta     \right],
\\
\\
\displaystyle
u_{tt}(0,t)= -\frac 1 a f'(t)  - \frac 2 a \sum_{k=1}^\infty \left[-\lambda_k e^{-\lambda_k t}f(0) 
+ e^{-\lambda_k t}f'(0) +  \int_0^t e^{-\lambda_k\eta}f''(t-\eta)d\eta     \right].
\end{array}
$$
Using $f(0)=f'(0)=0$, (2.5), and $\zeta(2) =  \frac {\pi^2} 6$, we have 
$$\begin{array}{c}
\displaystyle
\vert u_{t}(0,t)\vert
\le
\frac{1}{a_L}\vert f(t)\vert + \frac {a_U} 3 \Vert f'\Vert_{L^\infty[0,T]},
\end{array}
\tag{2.8}
$$
$$\begin{array}{c}
\displaystyle
\vert u_{tt}(0,t)\vert
\le
\frac{1}{a_L}\vert f'(t)\vert +\frac {a_U} 3 \Vert f'\Vert_{L^\infty[0,T]}.
\end{array}
\tag{2.9}
$$
Taking the supremum for (2.6),(2.8), and (2.9) for $t\in[0,T]$, we have 
$$\begin{array}{c}
\displaystyle
  \Vert u(0,\cdot) \Vert_{W^{2,\infty}[0,T]} \le \max\left(
\frac{\Vert f \Vert_{L^1[0,T]}}{a_L} + \frac{ a_U \Vert f \Vert_{L^\infty[0,T]}}  3, \right.\\
\\
\displaystyle
\left.
\frac{\Vert f \Vert_{L^\infty[0,T]}}{a_L} + \frac{ a_U \Vert f' \Vert_{L^\infty[0,T]}}  3,
\frac{\Vert f' \Vert_{L^\infty[0,T]}}{a_L} + \frac{ a_U \Vert f'' \Vert_{L^\infty[0,T]}}  3
\right).
\end{array}
$$
for all $0\le t\le T$. From this inequality, (2.3) and (2.4) follows. 
$\square$\\

\proclaim{\noindent Lemma 2.3}
If $0 < a \le a_U$ and $f\in\ L^\infty[0,\,T]$, then we have
$$
\Vert u(\,\cdot\,,T)\Vert_{L^1[0,\,a]}
\le C_T ||f||_{L^\infty[0,T]},
\tag {2.10}
$$
where
$$
C_T = \left( \frac 1 3  + \frac 2{3\pi} \right) a_U^2 + T.
$$
Further if $f(t)=t^r$, then
$$
\Vert u(\,\cdot\,,T)\Vert_{L^1[0,\,a]} \le C_{T,r},
\tag {2.11}
$$
where
$$
C_{T,r} =   \left( \frac 1 3  + \frac 2{3\pi} \right) T^r a_U^2 + \frac {T^{r+1}}{r+1}.
$$
\endproclaim

{\bf Proof of Lemma 2.3}\\
From  Lemma 2.1, we obtain
$$\begin{array}{c}
\displaystyle
u(x,T)=-\frac{a}3 f(T) \left(\frac 3{2a^2} x^2 - \frac 3 a x\right)
-\frac{1}{a}
\int_0^T f(t)dt\\
\\
\displaystyle
-\frac 2 a\sum_{k=1}^{\infty}\cos\left(\frac{k\pi}a x\right)\int_0^T e^{-\left(\frac{k\pi}{a}\right)^2(T-s)}f(s)ds.
\end{array}
$$
Therfore, by using (2.5) and $\zeta(2) = \frac {\pi^2} 6$, we have
$$
\begin{array}{ccl}
\displaystyle
\Vert u(\cdot,T)\Vert_{L^1[0,a]}
&\le& 
\displaystyle
\frac{a^2}3 \vert f(T) \vert 
+\Vert f\Vert_{L^1[0,T]}
+\frac{4}{\pi}\sum_{k=1}^{\infty}
\int_0^T e^{-\left(\frac{k\pi}{a}\right)^2(T-s)}\vert f(s)\vert ds
\\  
\\ 
\displaystyle
&\le& 
\displaystyle
\frac{a^2}3 \vert f(T) \vert
+\Vert f\Vert_{L^1[0,T]}
+\frac{2a^2}{3\pi}  \Vert f\Vert_{L^\infty[0,\,T]}
\\
\\
\displaystyle
&\le& 
 C_T \Vert f\Vert_{L^\infty[0,\,T]}.
\end{array}
$$
If $f(t)=t^r$, using $f(T) = \Vert f\Vert_{L^\infty[0,T]} = T^r$ and 
$\Vert f\Vert_{L^1[0,T]} = \frac{T^{r+1}}{r+1}$, we have the upper bound (2.11).
$\square$\\

For Lemma 2.4 and 2.6. let us define 
$$ b_j(r) = (-1)^{r-j} \frac{r!}{j!}, $$
for $r=0,1,2,\cdots$ and $0\le j\le r$.\\

\proclaim{\noindent Lemma 2.4}
If $f(t)=t^r,r=0,1,2,\cdots$, we have
\[ 
\begin{array}{ccl}
\displaystyle
u(0,t) &=&
\displaystyle
-\frac 1 {(r+1)a} t^{r+1} - \frac a 3 t^r - \frac 2 a  \sum_{j=0}^{r-1} b_j(r) t^j 
\left( \frac a\pi\right)^{2r+2-2j}\zeta(2r+2-2j)
\\
\\
\displaystyle
&&
\displaystyle
 - \frac {2(-1)^{r-1}r!} a \sum_{k=1}^\infty \frac 1{\lambda_k^{r+1}} e^{-\lambda_k t}. 
\end{array}
\]
\endproclaim

{\bf Proof of Lemma 2.4}\\
For $r\ge 0$, by induction argument, we have
\[
\int_0^t e^{\lambda (s-t)} s^r ds =  \sum_{j=0}^r \frac {b_j(r)}{\lambda^{r+1-j}} t^j - \frac {b_0(r)}{\lambda^{r+1}}e^{-\lambda t}.
\] 
From this formula and using
$$
\sum_{k=1}^\infty \frac 1{\lambda_k^{r+1-j}} = \left( \frac a\pi\right)^{2r+2-2j}\zeta(2r+2-2j),
$$
we obtain the lemma.
$\square$\\

For example,  for $r=0,1,2$ we have
\[
u(0,t) =
\left\{
\begin{array}{lc}
\displaystyle
-\frac 1 a t  - \frac a 3 
+ \frac 2 a \sum\limits_{k=1}^{\infty} \frac 1 {\lambda_k} e^{-\lambda_k t} 
& \text{if}\, f(t)=1;
\\
\\
\displaystyle
-\frac 1{2a} t^2 - \frac a 3 t  + \frac{a^3}{45}  
-\frac 2 a \sum\limits_{k=1}^\infty \frac 1 {\lambda_k^2} e^{-\lambda_k t} 
& \text{if}\, f(t)=t;
\\
\\
\displaystyle
- \frac{1}{3a} t^3  -\frac a 3 t^2 +\frac{2a^3}{45} t -\frac{4a^5}{945} 
+\frac 4 a \sum\limits_{k=1}^\infty \frac 1  {\lambda_k^3} e^{-\lambda_k t}    
& \text{if}\, f(t)=t^2.
\end{array}
\right.
\tag{2.12}
\]

\proclaim{\noindent Remark 2.5}
Here, we remark on the complexity of the correspondence  $a \longmapsto u(0,\cdot) \vert_{[0,T]}$. These examples suggest $u(0,\,\cdot\,)\vert_{[0,\,T]}$ for {\it general} $f$ contains information about $a$ that is quite complicated.
For example, when $f(t)=t^r$, $u(0,t)= O\left(\frac 1{a^{2r+1}}\right)$ by Lemma 2.4, resulting in large perturbation of Dirichlet data $\Delta u(0,t)$ from even in small negative perturbation of $\Delta a$, especially for small $a$ and large $r$ .
However, the enclosure method is not affected by the complexity and nonlinearity of the correspondence $a\longmapsto u(0,\,\cdot\,)\vert_{[0,\,T]}$ and yields $a$ explicitly, in particular,  with an explicit error estimate.
\endproclaim

Let us define the truncated approximation $u^N(t)$ of $u(0,t)$ in (2.2) as follows:
\[ 
u^N(t) = - \frac 1 a\int_0^t f(t) dt - 
\frac 2 a \sum_{k=1}^N  \int_0^t e^{-\lambda_k (t-s)} f(s) ds.
\tag{2.13}
\]
If $f\in W^{1,\infty}[0,T]$ and $f(0)=0$, let us change (2.2) and (2,13) as follows, by using integration by parts and 
$\sum_{k=1}^\infty \frac 1{\lambda_k}  = \frac{a^2}6 $:
\[
u(0,t) =
-\frac a 3 f(t) - \frac 1 a\int_0^t f(t) dt + \frac 2 a \sum_{k=1}^\infty \frac 1{\lambda_k} \int_0^t e^{-\lambda_k (t-s)} f'(s) ds;
\tag{2.14}
\]
\[ 
u^N(t) = -\frac a 3 f(t) - \frac 1 a\int_0^t f(t) dt + \frac 2 a \sum_{k=1}^N \frac 1{\lambda_k} \int_0^t e^{-\lambda_k (t-s)} f'(s) ds.
\tag{2.15}
\]

Then, the error between $u(0,t)$ and $u^N(t)$ is bounded by the following lemma:
\proclaim{\noindent Lemma 2.6}
If $f\in L^\infty[0,T]$, then
\[
\Vert u(0,\cdot)-u_N(\cdot)\Vert_{L^\infty[0,T]} \le \frac {2a}{\pi^2 N} ||f||_{L^\infty [0,T]}.
\]
Furthermore, if $f\in W^{1,\infty}[0,T]$ and $f(0)=0$, then 
\[
\Vert u(0,\cdot)-u_N(\cdot)\Vert_{L^\infty[0,T]} \le \frac {2a^3}{\pi^4 N^3} ||f||_{W^{1,\infty} [0,T]}.
\]
\endproclaim 

{\bf Proof of Lemma 2.6}\\
If $f\in L^\infty[0,T]$, using (2.2) and (2.13), then for all $t\in[0,T]$
\begin{eqnarray*}
&&|u(0,t)-u_N(t)| =\left|\frac 2 a \sum_{k=N+1}^{\infty} \int_0^t e^{-\lambda_k (t-s)} f(s) ds\right|
\\
&\le& \frac 2 a ||f||_{L^\infty [0,T]}\sum_{k=N+1}^{\infty} \int_0^t e^{-\lambda_k (t-s)} ds
\le \frac 2 a ||f||_{L^\infty [0,T]} \sum_{k=N+1}^{\infty} \frac 1  {\lambda_k}
\\
&\le& \frac {2a}{\pi^2} ||f||_{L^\infty [0,T]} \int_N^\infty  \frac 1 {x^2} dx
\le  \frac {2a}{\pi^2} ||f||_{L^\infty [0,T]}  \frac 1 N.
\end{eqnarray*}
If $f\in  W^{1,\infty}[0,T]$, using (2.14) and (2.15), then for all $t\in[0,T]$
\begin{eqnarray*}
&&|u(0,t)-u_N(t)| =\left|\frac 2a \sum_{k=N+1}^{\infty} \frac 1{\lambda_k} \int_0^t e^{-\lambda_k (t-s)} f'(s) ds\right|
\\
&\le&\frac 2a||f||_{W^{1,\infty} [0,T]}\sum_{k=N+1}^{\infty}  \frac 1{\lambda_k^2}
\le\frac{2a^3}{\pi^4}||f||_{W^{1,\infty} [0,T]}\int_N^\infty  \frac 1 {x^4} dx
\\
&\le&\frac{2a^3}{\pi^4}||f||_{W^{1,\infty} [0,T]} \frac 1{N^3}.
\end{eqnarray*} 
$\square$

For example, if $f(t)=t^2 e^{-\nu t},\nu>0$, then
\[
\begin{array}{l}
\displaystyle
    u_N(t) = -\frac a 3 t^2 e^{-\nu t} 
      +\frac 1 a \left( \frac {t^2} \nu + \frac {2t} {\nu^2} + \frac 2 {\nu^3}\right) e^{-\nu t} - \frac 2{a\nu^3} 
\\
\\
\displaystyle
     +\frac 2 a \sum_{k=1}^N \left[\left(  \frac {t^2}{\lambda_k}  
      -  \frac {t^2}{\lambda_k-\nu} + \frac {2t}{(\lambda_k-\nu)^2}  - \frac 2{(\lambda_k-\nu)^3}\right) e^{-\nu t}
      + \frac {2e^{-\lambda_k t}} {(\lambda_k-\nu)^3} e^{-\lambda_k t} \right].
\end{array}
\tag{2.16}
\]
Moreover, if $f(t)=t^r,r=0,1,2,\cdots$, we have
\[ 
\begin{array}{ccl}
u_N(t) &=& 
\displaystyle
-\frac {t^{r+1}} {(r+1)a} - \frac{at^r} 3 
- \frac {2(-1)^{r-1}r!} a \sum_{k=1}^N \frac {e^{-\lambda_k t}}{\lambda_k^{r+1}}
\\
\\
&&
\displaystyle
- \frac 2 a  \sum_{j=0}^{r-1} b_j(r) t^j 
\left( \frac a\pi\right)^{2r+2-2j}\zeta(2r+2-2j) . 
\end{array}
\tag{2.17}
\]
Then, we have the following truncation error:\\

\proclaim{\noindent Lemma 2.7}
For $f(t)=t^r,r=0,1,2,\cdots$, we have
$$
    |u(0,t)-u_N(t)| \le \frac {2(r!)a^{2r+1}}{(2r+1)\pi^{2r+2}} N^{-2r-1}. 
$$
If $t>\frac{a^2 \log(2)}{\pi^2(2N+3)}$, then
$$
   |u(0,t)-u_N(t)| \le \frac {4(r!)a^{2r+1}}{\pi^{2r+2}}\frac{2^{-\frac{(N+1)^2}{2N+3}}}{(N+1)^{r+1}}.
$$ 
\endproclaim

That is, the truncation error is of the order $O(N^{-2r-1})$ with a hyperconvergence of the order $O(2^{-N}N^{-3})$ for $t>\frac{a^2 \log(2)}{\pi^2(2N+3)}$.\\

{\bf Proof of Lemma 2.7}\\
Since
\[|u(0,t)-u_N(t)| =\frac {2 (r!)}{a}\sum_{k=N+1}^{\infty} \frac 1{\lambda_k^{r+1}} e^{-\lambda_k t}.
\]
and $e^{-\lambda_k t} \le 1$, we have
\begin{eqnarray*}
    |u(0,t)-u_N(t)|       & \le & \frac {2(r!)}{a} \left(\frac {a^2}{\pi^2}\right)^{r+1} \sum_{k=N+1}^{\infty} \frac 1 {k^{2r+2}}\\
       &\le& \frac {2(r!)}{a} \left(\frac {a^2}{\pi^2}\right)^{r+1} \int_N^\infty x^{-2M-2} dx \\
      &=&  \frac {2(r!)}{(2r+1)a} \left(\frac {a^2}{\pi^2}\right)^{r+1}  N^{-2r-1}. 
\end{eqnarray*}
Further, $t>\frac{a^2 \log(2)}{\pi^2(2N+3)}$ implies $s=e^{-\lambda_{N+2} t}/e^{-\lambda_{N+1}t}<\frac 1 2$.
For this $t$, we have 
$$
\sum_{k=N+1}^{\infty} e^{-\lambda_k t}
= \sum_{k=0}^{\infty} e^{-\lambda_{N+1} t} s^k \le 2 e^{-\lambda_{N+1} t}
$$
and 
\[|u(0,t)-u_N(t)| \le \frac {2 (r!)}{a}\frac 1{\lambda_{N+1}^{r+1}} 2 e^{-\lambda_{N+1} t}
\le  \frac {4(r!)a^{2r+1}}{\pi^{2r+2}}\frac{2^{-\frac{(N+1)^2}{2N+3}}}{(N+1)^{r+1}}.
\]
This proves Lemma 2.7.
$\square$\\

\section{Proof of Theorems 1.1 and 1.2}

\subsection{Proof of Theorem 1.1}
For the proof, we introduce $a^\infty(\tau)$ and divide the left side as follows:
$$
    \vert a- a(\tau)\vert \le \vert a-a^\infty(\tau)\vert + \vert a^\infty(\tau) - a(\tau) \vert.
$$
Let us first introduce $\hat u^\infty(\tau), I^\infty(\tau),$ and $a^\infty(\tau)$. It is not to difficult to show that $\hat u (\cdot,\tau)$ is the unique solution of the boundary value problem
$$
\left\{
\begin{array}{l}
w''-\tau^2 w=e^{-\tau^2 T}u(\cdot,T)\,\,\text{in}\,]0,\,a[,\\
\\
\displaystyle
w'(0)=\hat{f}(\tau),\,\,w'(a)=0.
\end{array}
\right.
\tag {3.1}
$$
Let $\hat u^{\infty}=\hat u^{\infty}(\cdot,\tau)$ be the unique solution of the boundary value problem
$$
\left\{
\begin{array}{c}
w''-\tau^2 w=0\,\,\text{in}\,]0,\,a[,\\
\\
\displaystyle
w'(0)=\hat{f}(\tau),\,\,w'(a)=0,
\end{array}
\right.
\tag {3.2}
$$
and  $\hat u^{\infty}$ has the explicit expression
$$\displaystyle
\hat u^{\infty}(x,\tau)
=-\frac{\hat{f}(\tau)}{\tau\sinh (a\tau)}
\cosh (\tau(x-a)).
\tag {3.3}
$$
Recalling that, from (1.3),  
$$\displaystyle
I(\tau)
=\tau \hat u(0,\tau)+\hat{f}(\tau),
$$
let us define 
$$\displaystyle
I^{\infty}(\tau)
=\tau \hat u^{\infty}(0,\tau)+\hat{f}(\tau),
\tag {3.4}
$$
$$\displaystyle
a^{\infty}(\tau)=\frac{1}{-2\tau}\log\left\vert\frac{I^{\infty}(\tau)}{\displaystyle
-2\hat{f}(\tau)}\right\vert.
\tag {3.5}
$$
Inserting (3.3) into (3.4), we obtain
$$\displaystyle
I^{\infty}(\tau)
=-2\hat{f}(\tau)\frac{e^{-2a\tau}}{1-e^{-2a\tau}}.
\tag {3.6}
$$
This, together with (3.5) yields
$$\displaystyle
e^{-2a^{\infty}(\tau)\tau}
=\frac{e^{-2a\tau}}{1-e^{-2a\tau}}
$$
and thus
$$\displaystyle
e^{2(a-a^{\infty}(\tau))\tau}
=\frac{1}{1-e^{-2a\tau}}.
$$
Taking the logarthm on both sides, we obtain
$$\displaystyle
a-a^{\infty}(\tau)
=\frac{-\log(1-e^{-2a\tau})}{2\tau}.
\tag {3.7}
$$
Here, we note that, for all $x\in\,]0,\,1[$,
$$\displaystyle
0<\log\frac{1}{1-x}<\frac{x}{1-x},
$$
and the function $-\log(1-e^{-2a\tau})/\tau, \tau>0$ decreases monotonically.
These results, together with (3.7), yield, for all $\tau\ge\tau_0$
$$\displaystyle
0<a-a^{\infty}(\tau)
\le
\frac{-\log(1-e^{-2a\tau_0})}{2\tau_0}
\le
\frac{-\log(1-e^{-2a_L\tau_0})}{2\tau_0},
\tag{3.8}
$$
where $\tau_0$ is an arbitrary positive number. 

\[
\]
Next we provide an upper estimate for $\vert a(\tau)-a^{\infty}(\tau)\vert$.
\noindent
Because
$$\displaystyle
I(\tau)=I^{\infty}(\tau)+\tau (\hat u(0,\tau)-\hat u^{\infty}(0,\tau)),
$$
by defining
$$\displaystyle
E(\tau)
=\frac{\tau(\hat u(0,\tau)-\hat u^{\infty}(0,\tau))}{I^{\infty}(\tau)},
$$
we have
$$
\left\{
\begin{array}{l}
\displaystyle
I(\tau) =I^{\infty}(\tau)(1+E(\tau)),\\
\\
\displaystyle
a(\tau) =a^{\infty}(\tau)
-\frac{1}{2\tau}\log\vert 1+E(\tau)\vert.
\end{array}
\right.
\tag {3.9}
$$
Using the method of variation of parameters,  the solution $\hat u $ of nonhomogeneous ordinary differential equation (3.1) with a Neumann boundary condition could be computed as :
$$\begin{array}{l}
  \hat u (x,\tau) = \hat u^\infty(x,\tau) - e^{-\tau^2 T}\times\\
\\
\displaystyle
\left(
\frac{\int_0^a u(\xi,T)\cosh (\tau\xi)d\xi}{\tau\sinh(a\tau)}\cosh (\tau(x-a))
+\frac{1}{\tau}\int_x^au(\xi, T)\sinh (\tau(x-\xi))d\xi
\right).
\end{array}
$$
From this, we have
$$\begin{array}{ccl}
&&\vert \hat u (x,\tau) - \hat u^\infty(x,\tau)\vert 
\\
\\
&\le& 
\displaystyle
e^{-\tau^2 T}
\left\vert
\frac{\int_0^a u(\xi,T)\cosh (\tau\xi)d\xi}{\tau\sinh(a\tau)}\cosh (\tau a)
-\frac{1}{\tau}\int_0^au(\xi, T)\sinh (\tau\xi))d\xi
\right\vert
\\
\\
\displaystyle
&=& 
\displaystyle
\frac{e^{-\tau^2 T}}{\tau\sinh (a\tau)}
\left\vert\int_0^a u(\xi,T)\cosh (\tau(\xi-a))d\xi\right\vert
\\
\\
&\le& 
\displaystyle
\frac{e^{-\tau^2 T}\cosh (a\tau)}{\tau\sinh (a\tau)}
\Vert u(\,\cdot\,,T)\Vert_{L^1(0,\,a)}.
\end{array}
\tag{3.10}
$$
Using (1.4),(2.3),(3.6),(3.9),(3.10) and $\coth(a\tau)\le e^{a\tau}$, we have
$$
\begin{array}{lll}
\vert E(\tau)\vert
&\le& 
\displaystyle
\frac{(e^{2a\tau}-1)e^{-\tau^2T}}{2\vert\hat{f}(\tau)\vert}\coth (a\tau) \Vert u(\,\cdot\,,T)\Vert_{L^1(0,\,a)}\\
\\
\displaystyle
&\le&  
\displaystyle
\frac{e^{-\tau^2T+3a\tau}\tau^\mu}{2 C_\mu} C_T \Vert f\Vert_{L^\infty[0,T]}.
\end{array}
\tag{3.11}
$$
For all $\eta\in\,]-1,\,1[$, it holds that
$$\displaystyle
\vert\log(1+\eta)\vert\le\frac{\vert \eta\vert}{1-\vert \eta\vert}.
\tag{3.12}
$$
By using (3.9), (3.10), and (3.12), we obtain
$$\begin{array}{c}
\displaystyle
\vert a(\tau)-a^{\infty}(\tau)\vert
\le
\frac{e^{-\tau^2T+3a\tau}\tau^{\mu-1}}{4 C_\mu (1-\epsilon)} C_T \Vert f\Vert_{L^\infty[0,T]}
\end{array}
\tag {3.13}
$$
provided that
$$\displaystyle
\frac{e^{-\tau^2T+3a\tau}\tau^\mu}{2 C_\mu} C_T \Vert f\Vert_{L^\infty[0,T]}
\le \epsilon < 1.
$$

Let us define
$$\displaystyle
\eta(\tau)=e^{-\tau^2T+3a\tau}\tau^\mu.
$$
Differentiating $\eta(\tau)$ with respect to $\tau$, we have
$$
\frac{\eta'(\tau)}{\eta(\tau)}
=\frac{-2T\tau^2  + 3a \tau + \mu}\tau. 
$$
Therefore, $\eta'(\tau)<0$ if
$$\displaystyle
\tau\ge\frac{3a_U + \sqrt{9a_U^2 + 8T\mu}}{4T}\ge\frac{3a + \sqrt{9a^2 + 8T\mu}}{4T}.
\tag {3.14}
$$
Therefore, if $\tau_0$ satisfies (3.14), using $\eta(\tau)$ as a decreasing function for $\tau\ge \tau_0$,
we have (3.13). From (3.8) and (3.13), we could prove the theorem. 
$\square$\\

\subsection{Proof of Theorem 1.2}
Applying (1.13) to $g(t)=e^{-\tau^2 t}u(0,t)$ and (2.3), we obtain
$$\displaystyle
\vert Q_{N_t}(\hat u (\tau))-\hat u(\tau)\vert
\le \frac{T^3 C_{\max}\Vert f\Vert_{W^{2,\infty}[0,T]}}{12N_t^2}(\tau^2+1)^2.
\tag {3.15}
$$
From (1.4),(3.6), and (3.9), and by using $\vert E(\tau) \vert < \epsilon$, we obtain
$$\displaystyle
\vert I(\tau)\vert
\ge
\frac{e^{-2ar}}{1-e^{-2a\tau}} 2 C_\mu \tau^{-\mu}(1-\epsilon)
 =  \frac{2 C_\mu (1-\epsilon)}{\tau^\mu ( e^{2a\tau}-1 )}.
\tag {3.16}
$$
If we define 
$$
Z_{N_t}(\tau) = \frac{\tau( Q_{N_t}(\hat u (\tau))-\hat u(\tau) )}{I(\tau)}, 
$$
from (1.3) and (1.14), we have
$$
  I_{N_t}(\tau) = I(\tau) (1 + Z_{N_t}(\tau)).
$$
By using (3.15) and (3.16), we have
$$\displaystyle
\vert Z_{N_t}(\tau)\vert
\le\frac{T^3C_{\max}\Vert f\Vert_{W^{2,\infty}[0,T]}}{24 C_\mu (1-\epsilon) } \frac{\tau^{1+\mu}e^{2a\tau}(\tau^2+1)^2}{N_t^2}.
\tag {3.17}
$$
Now, using $N_t\ge N_t^{\delta}(\tau)$, it follows from (3.17) that,
for all $\tau\ge\tau_0$
$$\begin{array}{c}
\displaystyle
\vert Z_{N_t}(\tau)\vert
\le
\frac{T^3C_{\max}\Vert f\Vert_{W^{2,\infty}[0,T]}}{24 C_\mu (1-\epsilon) }\left(\frac 1{\tau^{\delta}}+\frac{1}{\tau^{2+\delta}}\right)^2.
\end{array}
\tag {3.18}
$$
Likewise (3.9), we have
$$
a_{N_t}(\tau) = a(\tau) - \frac 1{2\tau} \log\vert 1 + Z_{N_t}(\tau)\vert.
\tag {3.19}
$$
From (3.12) and (3.19), the inequality (1.16) is derived.
$\square$\\

\section{Numerical test}

In this section, the computation of the trusted frequency region for the enclosure method is presented, theoretically 
in Section 4.1 and numerically in Section 4.2 and 4.3.

\subsection{Trusted frequency region for $f(t)=t^2$ : Theoretical computation}
For $f(t)=t^2$, let us choose parameters $\tau_0, T, a_U, a_L, \delta $ satisfying (1.10),(1.11), and (1.15). We have chosen  
\[
\tau_0 = 3,\; T = 5,\; a_U=a_L=a=1, \; \delta=5.
\]
In this subsection, it is checked for these parameters to satisfy (1.10),(1.11), and (1.15). First, by (1.6), (2.4), and (2.11), we have:
\[
\left\{
\begin{array}{l}
\displaystyle
\mu=6; \\
\\
\displaystyle
C_\mu=
2(1-\frac {2.5}{e})\,\,\text{if}\,\tau_0 \ge \frac 1{\sqrt{ T}};\\
\\
\displaystyle
C_{T,2}=
\left( \frac 1 3 + \frac 2{3\pi}   \right) T^r a_U^2+\frac{T^{r+1}}{r+1};\\
\\
\displaystyle
C_{max,2}=
\frac{T^{r+1}}{a_L} + \frac{ T^r a_U} 3\,\,\text{if}\,T\ge 3.
\end{array}
\right.
\]
By Lemmas 2.2 and 2.3, $C_{T,2}$ and $C_{max,2}$ replace $C_T \Vert f \Vert_{L^\infty[0,T]}$ and $C_{max}\Vert f\Vert_{W^{2,\infty}[0,T]}$, respectively.

In Figures 1(a), 1(b), and 1(c), the validity of (1,10), (1.11), and (1.15) is checked by defining 
\[
\left\{
\begin{array}{l}
\displaystyle
F(T)=\frac{3 a_U}{4T} \left(1+\sqrt{1+\frac{8T\mu}{9 a_U^2}}\,\right),\\
\\
\displaystyle
G(\tau)=\frac{e^{-T\tau_0^2+3a_U\tau_0}\tau_0^\mu}{2C_\mu} C_T \Vert f\Vert_{L^\infty[0,T]},\\
\\
\displaystyle
H(\tau)=
\frac{T^3C_{\max}\Vert f\Vert_{W^{2,\infty}[0,T]}}{24 C_\mu (1-\epsilon)}\left(\frac{1}{\tau_0^{\delta}}+\frac{1}{\tau_0^{2+\delta}}\right)^2.
\end{array}
\right.
\]

$F(T)$ is plotted in Figure 1(a) and from this, we have
\[  F(T)= F(5)\le 1 < \tau_0 \]
and (1.10) is satisfied. The red horizontal line represents $1$.

In Figure 1(b), $G(\tau)$ is plotted, and by simple computation, we obtain
\[  G(\tau_0)=G(3)=\epsilon = 2.9114\times 10^{-11} \]
and  (1.11) is satisfied for this $\epsilon$. The red horizontal line represents $\epsilon$.

$H(\tau)$ is plotted in Figure 1(c). The red horizontal
line represents $\eta=H(3)=0.0904$.  For this $\eta$, (1.15) is satisfied.

$N_t^\delta(\tau)$ is plotted in Figure 1(d). By computation, we have $N_t^\delta(\tau_0)=2054266$, represented by the lower red horizontal line in the figure. If we choose $N_t = 10^10$, for example, we could choose $\tau_{max} = 5$. We could verify that $10^{10} > N_t^\delta(\tau_{max})=N_t^\delta(5)$ by examing the upper red horizontal line in the figure.


By using Theorem 1.3, we obtain the error bound in the trusted frequency region $[\tau_0,\tau_{max}]=[3,5]$ as follows:
\[
|a-a_{N_t}(\tau)|\le -\frac{\log(1-e^{-2a\tau_0})}{2\tau_0}+ \frac\epsilon{2\tau_0(1-\epsilon)} 
   +\frac\eta{2\tau_0(1-\eta)} \le 0.017.
\]

\begin{figure}
\begin{minipage}[t]{6cm}
\centerline{\epsfig{file=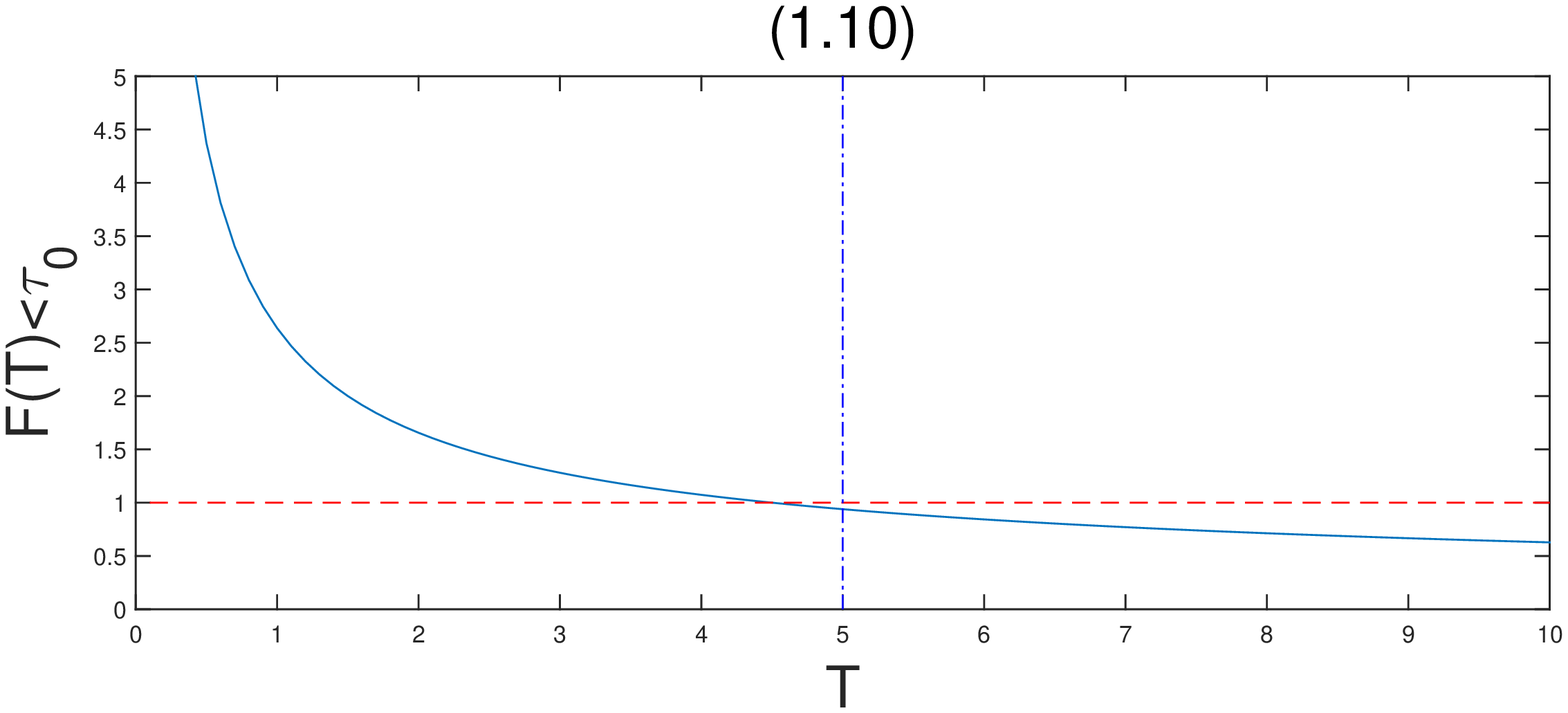,height=6cm,width=6cm,clip=}}
\begin{center}(a)\end{center}
\end{minipage}
\begin{minipage}[t]{6cm}
\centerline{\epsfig{file=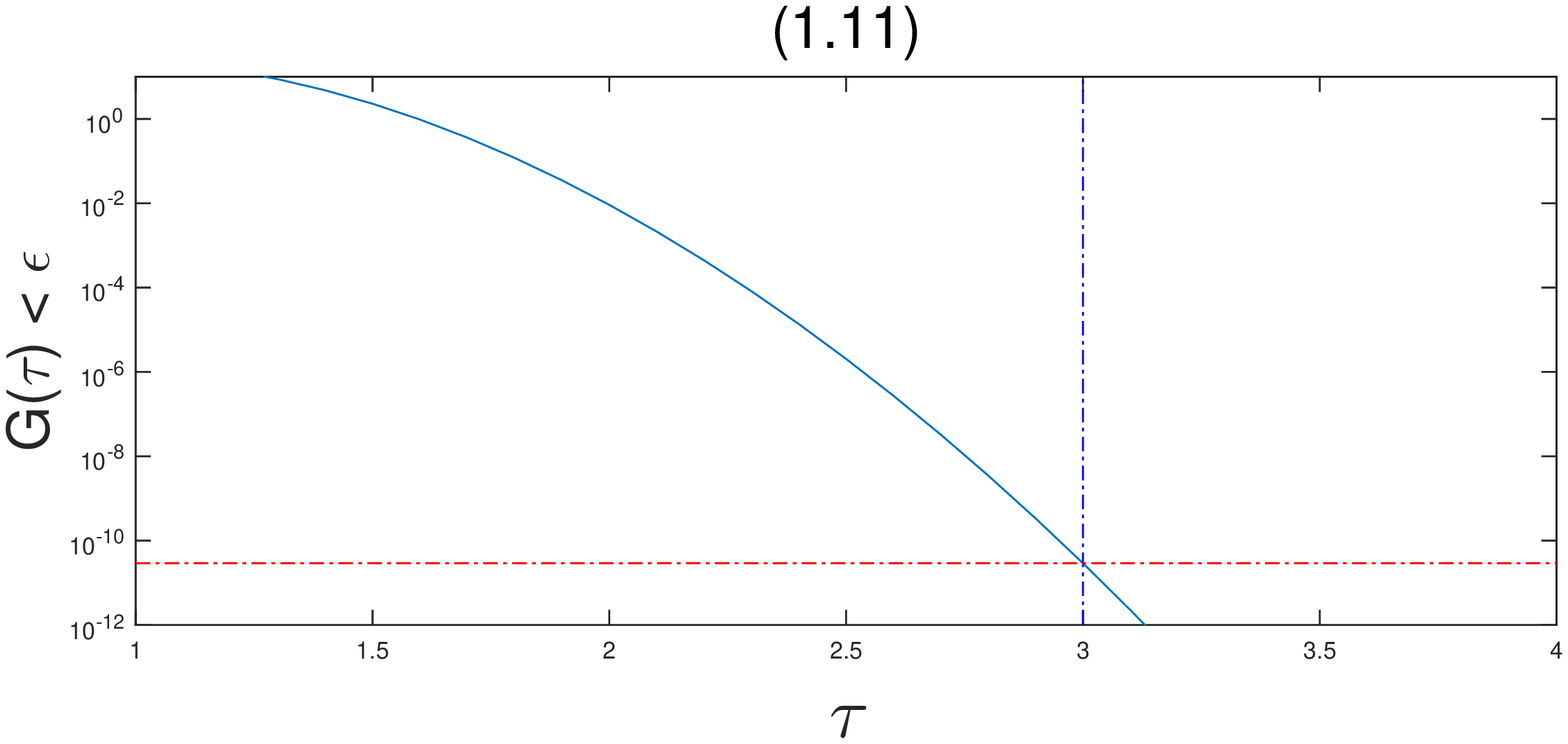,height=6cm,width=6cm,clip=}}
\begin{center}(b)\end{center}
\end{minipage}
\\
\\
\begin{minipage}[t]{6cm}
\centerline{\epsfig{file=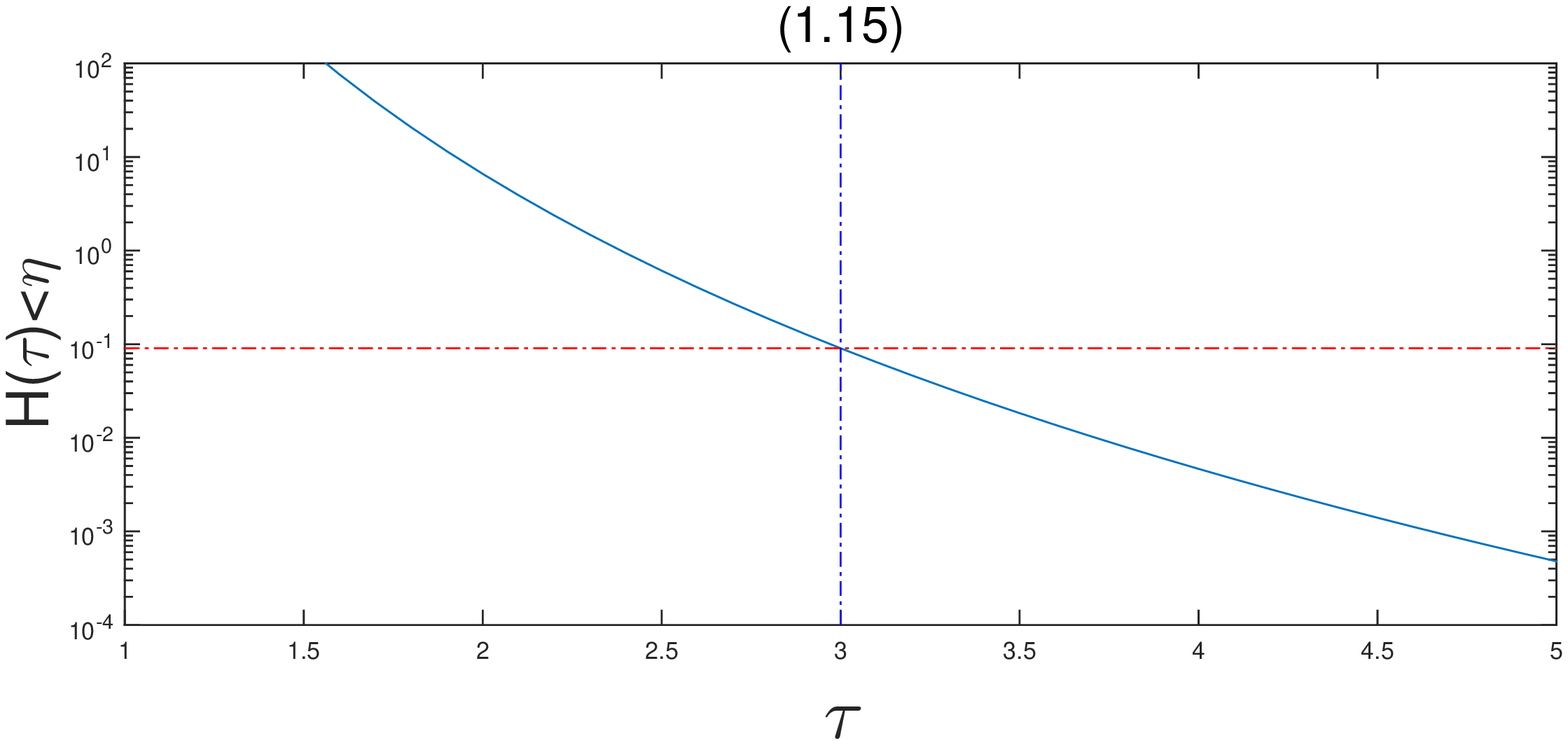,height=6cm,width=6cm,clip=}}
\begin{center}(c)\end{center}
\end{minipage}
\begin{minipage}[t]{6cm}
\centerline{\epsfig{file=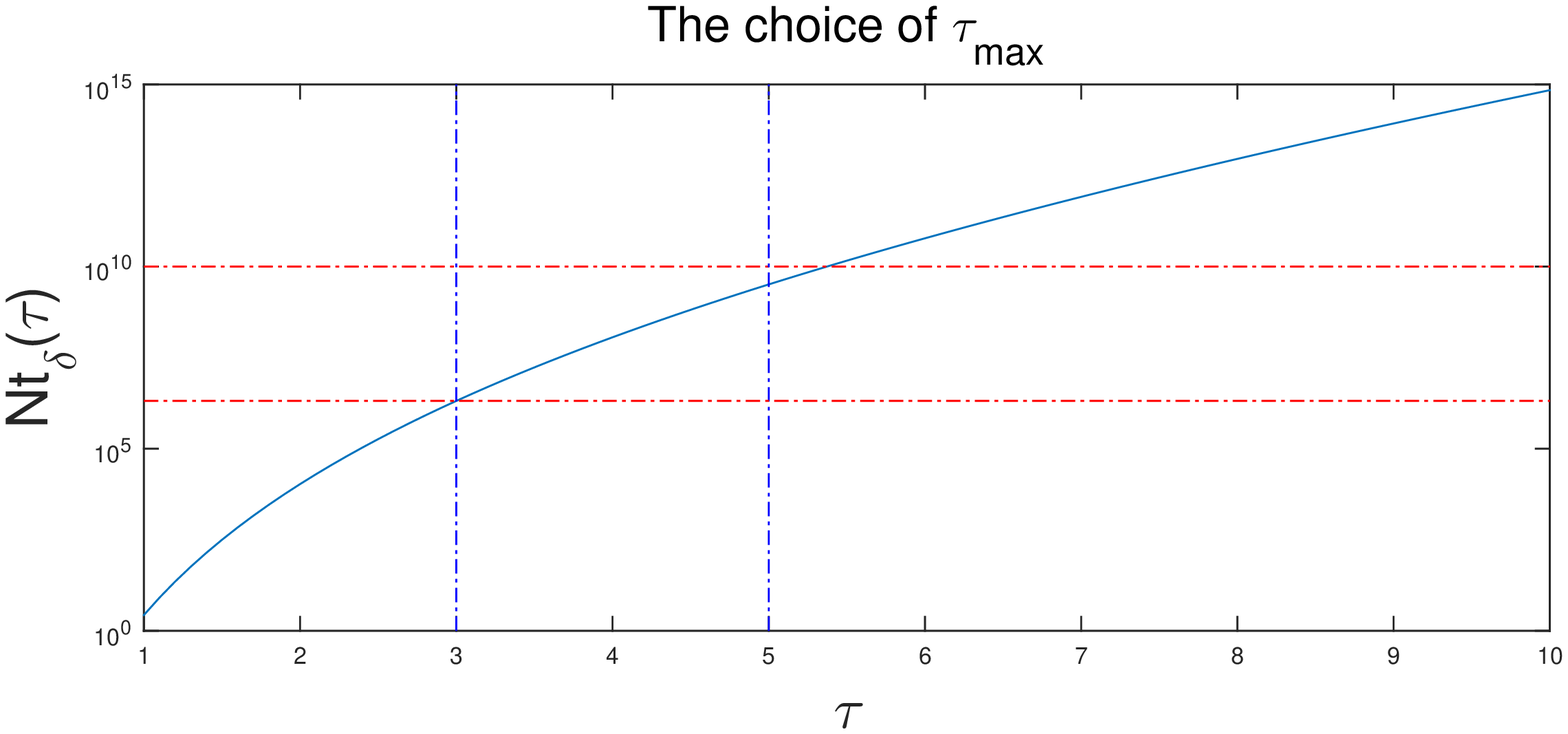,height=6cm,width=6cm,clip=}}
\begin{center}(d)\end{center}
\end{minipage}
\caption{ Computing the trusted frequency region $[\tau_0,\tau_{max}]$ for $f(t)=t^2$: (a) $F(T)$  (b) $G(\tau)$ (c) $H(\tau)$ (d) $N_t^\delta(\tau)$.  }
\end{figure}

\subsection{Trusted frequency region for $f(t)=t^r, r=0,1,2$ : Numerical computation}
The plots presented In Figure 1 enabled us to theoretically investigate the trusted frequency region with error bound 0.017. In this subsection, we numerically investigate the trusted frequency region with an error less than 0.01, (0.1 when $f(t)=1,t$),  in Sections 4.2 and 4.3. $T=5$ is fixed as in the previous section. The frequency is specified in increments 0.5 starting from 1, i.e. $1,1.5,2, \cdots $.

At first, in Figure 2(a), we fixed $a=1,N=10^3,f(t)=t^2$ and computed the trusted frequency region with an error bound of 0.01 
depending on $N_t=10^3,10^4,10^5,10^6$. The regions are as follows: 
\begin{itemize} 
\item{$N_t=10^3:  \;\;  [2.0, 5.0]$}
\item{$N_t=10^4:  \;\;  [2.0,  8.0]$}
\item{$N_t=10^5:  \;\;  [2.0, 11.0]$}
\item{$N_t=10^6:  \;\;  [2.0, 15.0]$}
\end{itemize}
This trusted frequency region is larger than the trusted frequency region that was theoretically determined in Section 4.1, even though $N_t$ is less than $10^10$ and the error bound is less than $0.017$.  The trusted frequency region for the same error bound becomes larger as $N_t$ increases.

Next, in Figure 2(b), we fixed $a=1,N=10^3, N_t=10^3$ and computed the trusted frequency region with an error bound of 0.1 depending on $f(t)=1,t,t^2$. The regions are as follows:
\begin{itemize}
\item{$f(t)=t^2: \;\;  [1.0, 6.0]$} 
\item{$f(t)=t :   \;\; [1.0, 2.0]$}
\item{$f(t)=1: \;\;$  None}
\end{itemize}
Here, the trusted frequency region becomes larger as $r$ increases. Note that $f(t)=t^2$ satisfies the assumption $f(0)=0, f'(0)=0$ in Theorem 1.2, whereas $f(t)=1$ and $f(t)=t$ does not. 

In Figure 2(c), we fixed $N=10^3, N_t=10^4, f(t)=t^2$ and computed the trusted frequency region with an error bound of 0.01  depending on $a=1,2,3,4$. The regions are as follows: 
\begin{itemize}
\item{$a=1: \;\;  [2.0, 8.0]$} 
\item{$a=2:   \;\; [2.0, 4.5]$}
\item{$a=3: \;\;  [2.5, 3.5]$}
\item{$a=4: \;\; [2.5, 2.5]$}
\end{itemize}
In this way, we established that the trusted frequency region becomes smaller as $a$ increases.

\begin{figure}
\begin{minipage}[t]{6cm}
\centerline{\epsfig{file=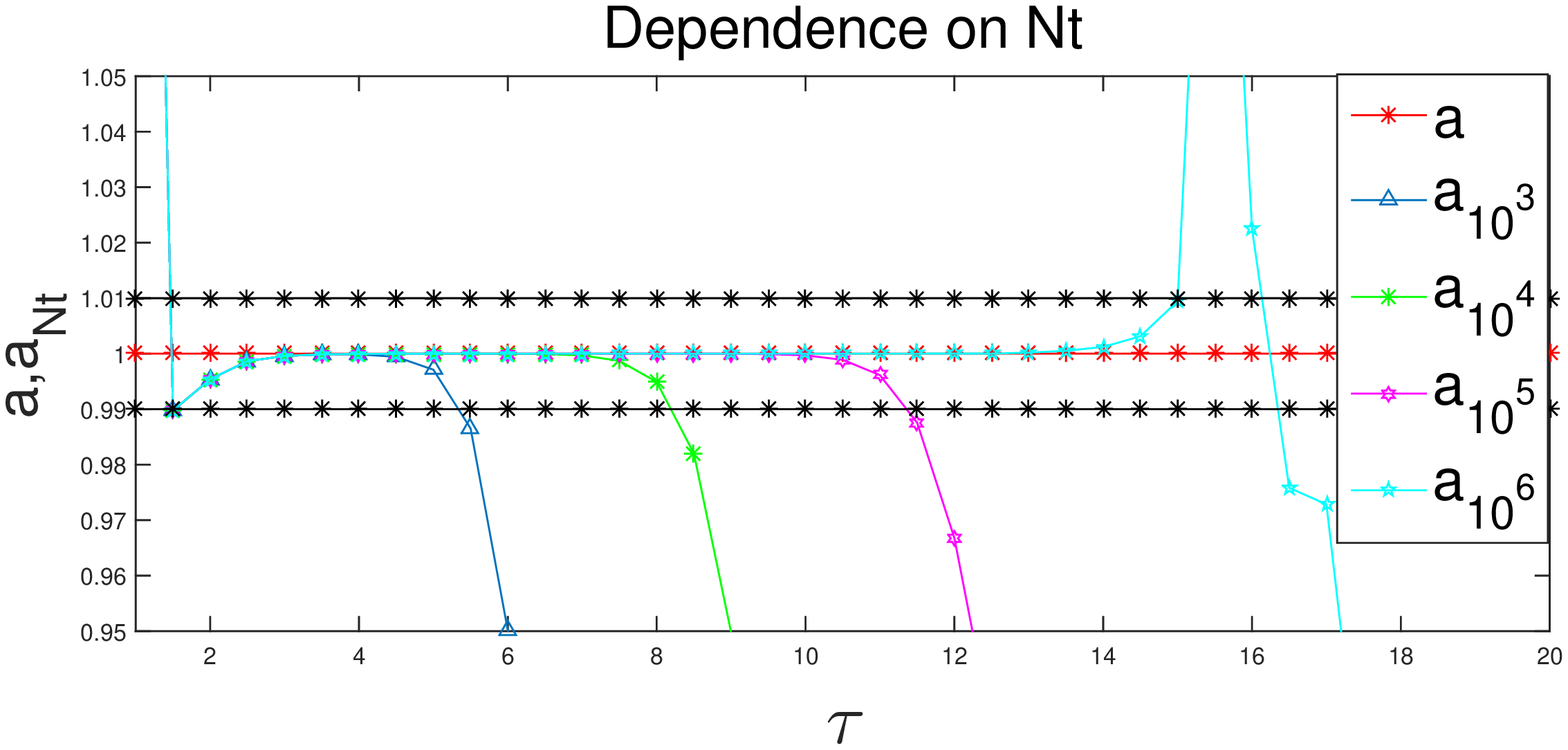,height=5cm,width=5cm,clip=}}
\begin{center}(a)\end{center}
\end{minipage}
\begin{minipage}[t]{6cm}
\centerline{\epsfig{file=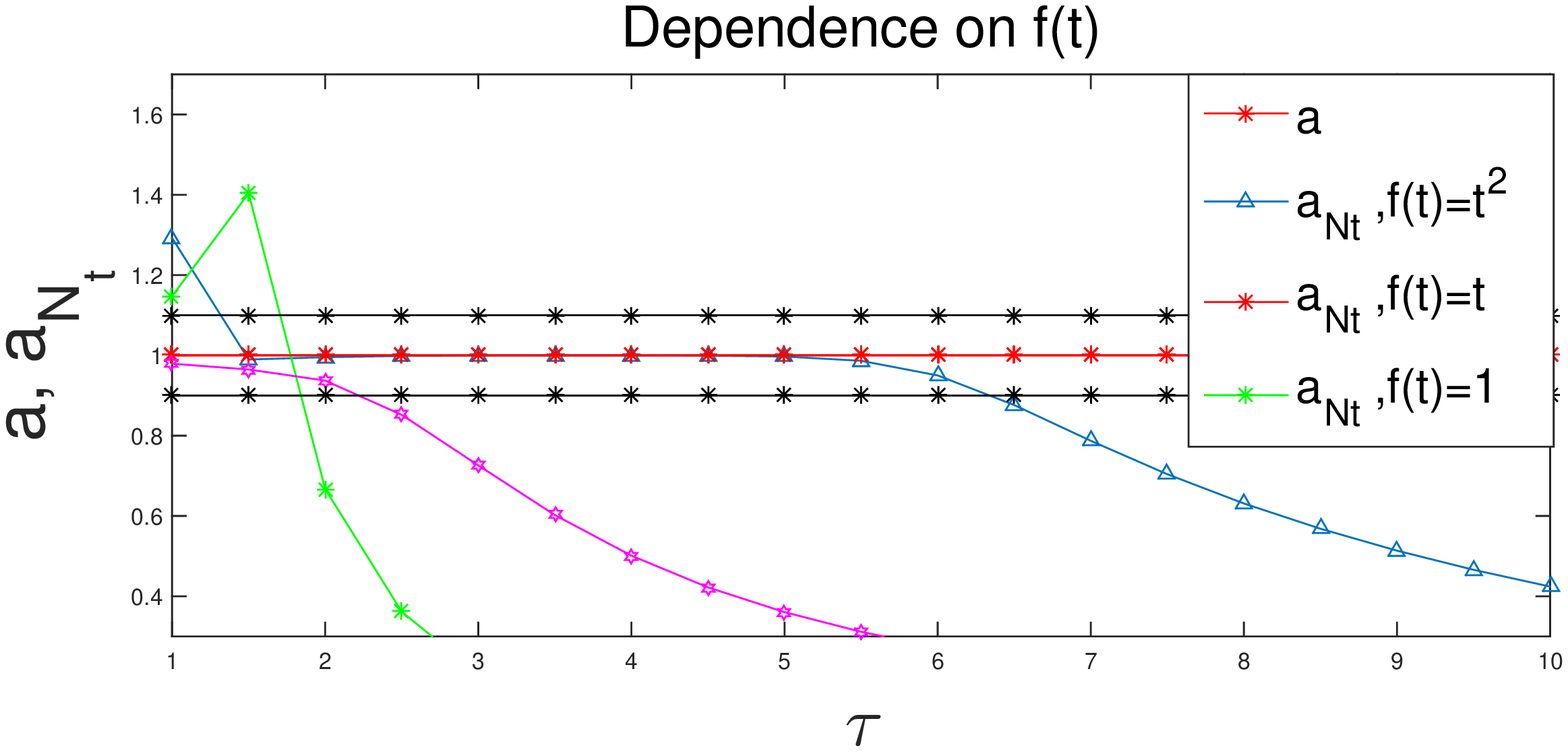,height=5cm,width=5cm,clip=}}
\begin{center}(b)\end{center}
\end{minipage}
\\
\begin{minipage}[t]{6cm}
\centerline{\epsfig{file=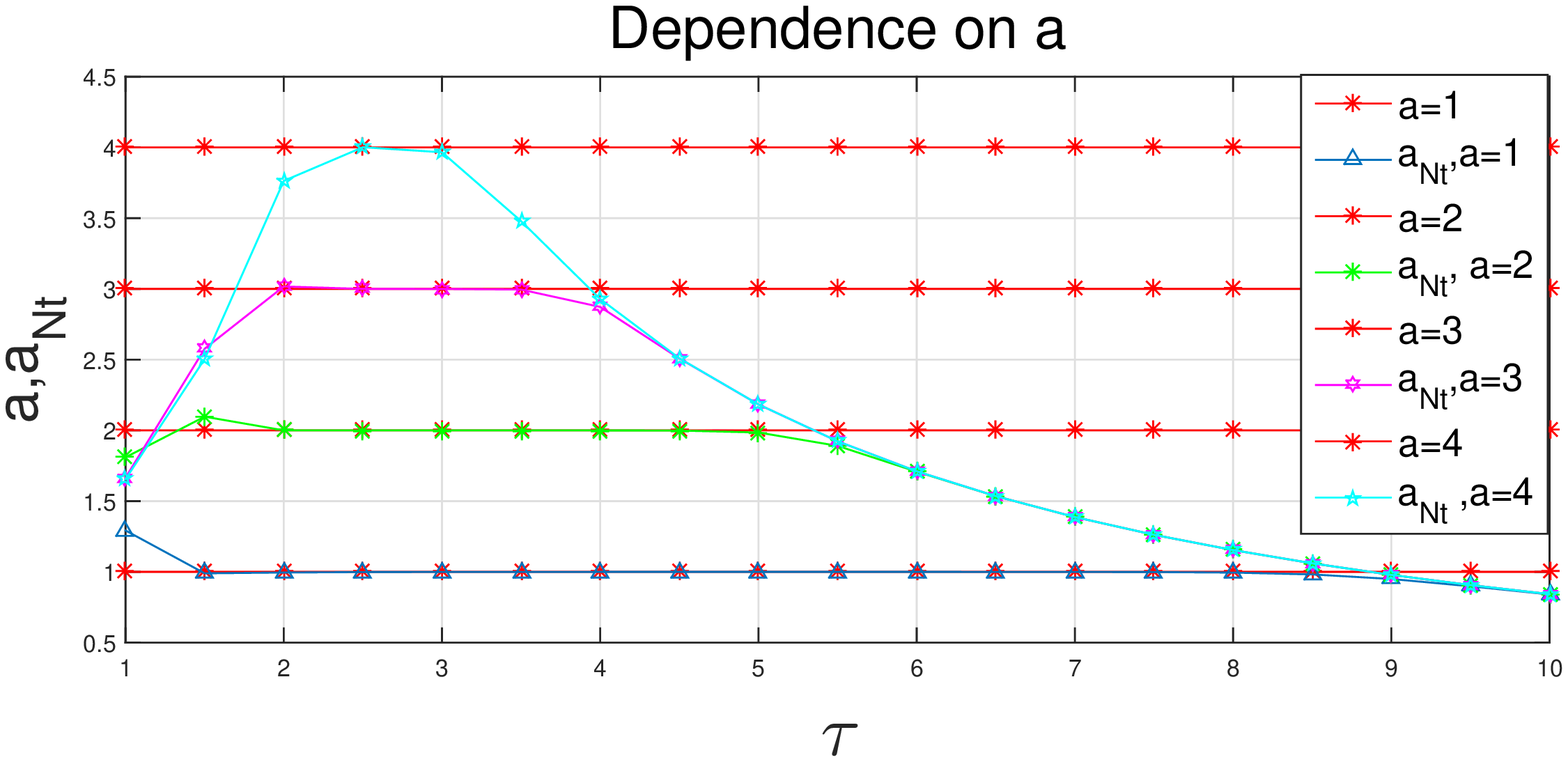,height=5cm,width=5cm,clip=}}
\begin{center}(c)\end{center}
\end{minipage}
\begin{minipage}[t]{6cm}
\centerline{}
\end{minipage}
\caption{ The approximation $a_{N_t}$ of true value $a$ with respect to $\tau$ depending on (a) the number of divisions for the time interval $N_t=10^3,10^4,10^5,10^6$ when $a=1,N=10^3,f(t)=t^2$ are fixed, (b) the right side function $f(t)=1,t,t^2$, when $a=1,N=10^3, N_t=10^3$ are fixed, and (c) the true value $a=1,2,3,4$, when $N=10^3, N_t=10^4, f(t)=t^2$ are fixed. }
\end{figure}

\subsection{Trusted frequency region for $f(t)=t^2e^{-2t}$ : Numerical computation} 
We also fixed $a=1,T=5,$ and $N=10^2$ in this subsection. We use (2.16) with $\nu=2$ for $u_N$, instead of (2.17) which is used when $f(t)=t^r$. Figure 3 (a) shows a plot of $e^2 f(t)$, in which $f(t)$ is normalised for the maximum value to be 1. 

In Figure 3(b), we fixed $a=1,N=10^3,f(t)=t^2 e^{2-25}$ and computed the trusted frequency region with an error bound of 0.01 depending on $N_t=10^3,10^4,10^5,10^6$. The regions are as follows: 
\begin{itemize} 
\item{$N_t=10^3:  \;\;  [2.0, 5.0]$}
\item{$N_t=10^4:  \;\;  [2.0, 8.0]$}
\item{$N_t=10^5:  \;\;  [2.0, 9.0]$}
\item{$N_t=10^6:  \;\;  [2.0, 9.0]$}
\end{itemize} 
The trusted frequency region becomes larger as $N_t$ increases as in the case $f(t)=t^2$ in Figure 2(a); however. the region 
is slightly smaller than that of $f(t)=t^2$, especially for $N_t=10^5, 10^6$. Moreover, the result for $N_t=10^6$ is slightly
better than that of $N_t=10^5$, although this was not discernable in the figure. The fact that the result for $f(t)=t^2$ is better than that for $f(t)=t^2 e^{-2t}$ comes from the approximation $u_N$ in (2.17), where the use of the Reimann-Zeta function improves the result compard to (2.16), where this function was not used. This inference can be verified by the order of convergence
$O(N^{-3})$ for $f(t)=t^2 e^{-2t}$ and  $O(N^{-5})$ (with hyperconvergence $O(2^{-N} N^{-3})$ for nonsmall $t$ )  for $f(t)=t^2$, which is shown in Lemma 2.6 and Lemma 2.7, respectively.

\begin{figure}
\begin{minipage}[t]{6cm}
\centerline{\epsfig{file=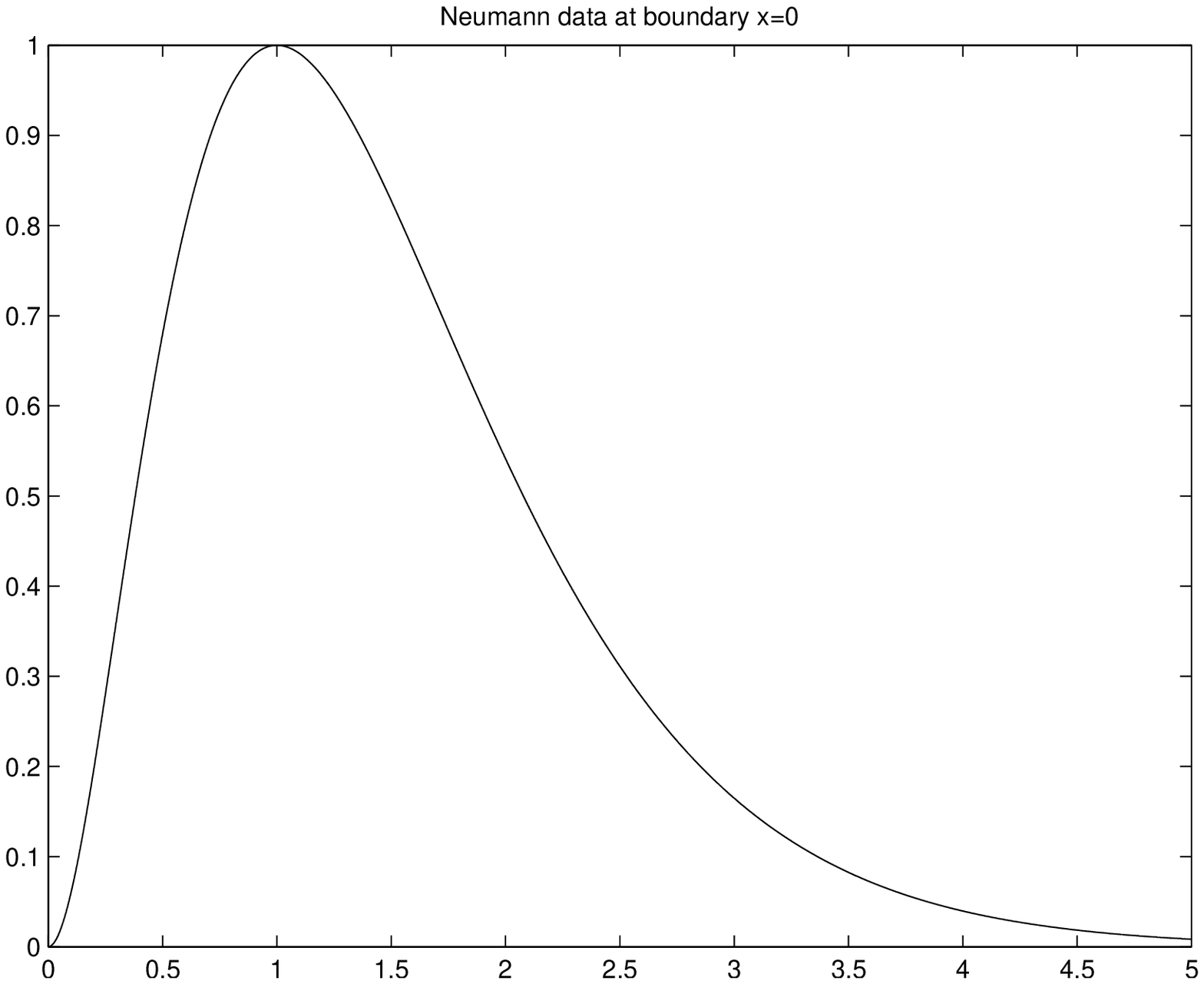,height=5cm,width=5cm,clip=}}
\begin{center}(a)\end{center}
\end{minipage}
\begin{minipage}[t]{6cm}
\centerline{\epsfig{file=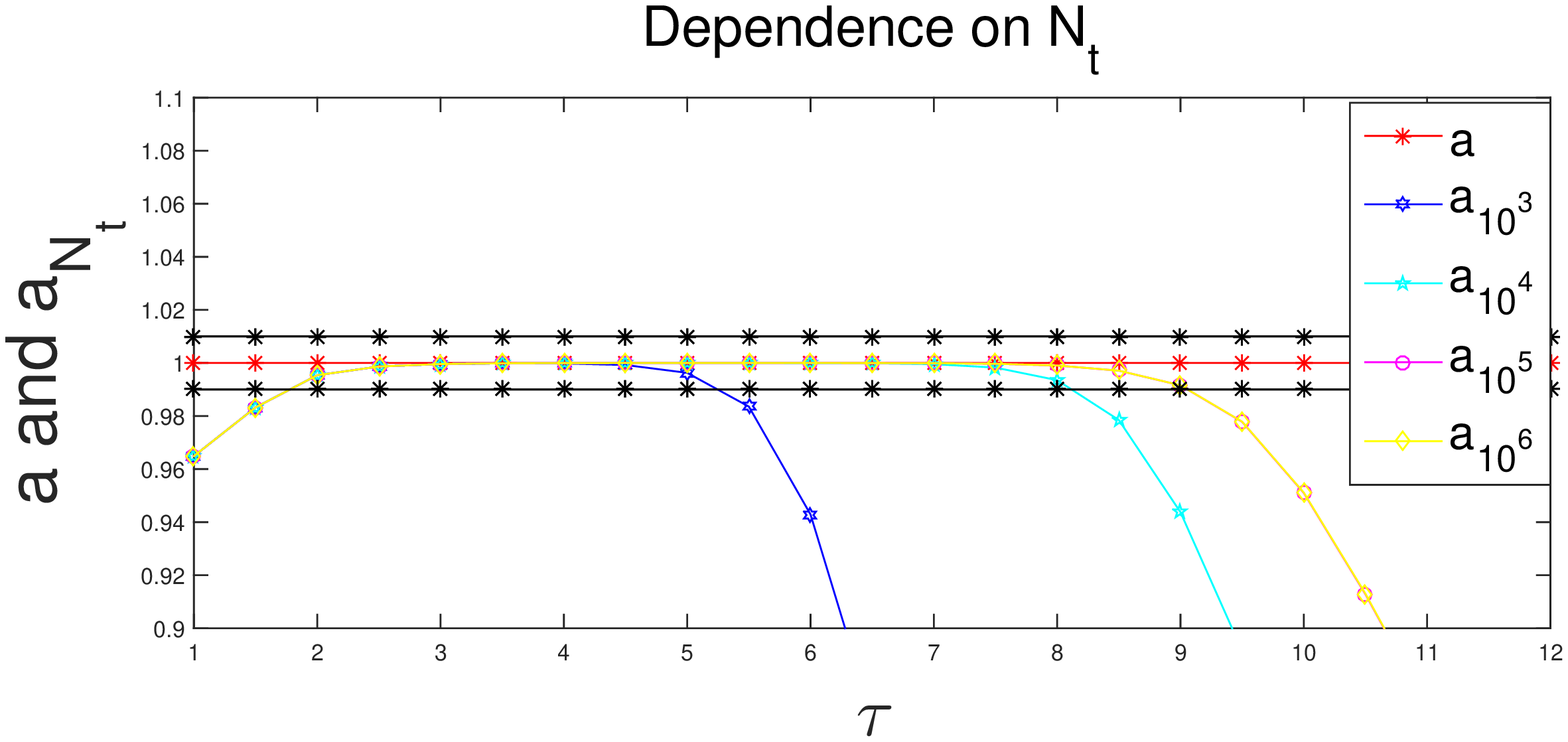,height=5cm,width=5cm,clip=}}
\begin{center}(b)\end{center}
\end{minipage}
\caption{(a) Neumann data $f(t)=  t^2 e^{2-2t}$ at $x=0$ and (b) the approximation of $a=1$, when $f(t)= t^2 e^{2-2t}, N=10^3, N_t=10^3,10^4,10^5,10^6$  }
\end{figure}

$$\quad$$

\centerline{{\bf Acknowledgements}}
This work 
was partially supported by Grant-in-Aid for
Scientific Research (C)(No. 25400155) of Japan  Society for the
Promotion of Science.
KK was also supported by the National Research Foundation of Korea (NRF) funded by
the Ministry of Education (NRF-2013R1A1A2010624).


\end{document}